\documentclass[12pt,english,a4paper]{amsart}
\usepackage[T1]{fontenc}
\usepackage[latin1]{inputenc}
\usepackage{geometry}
\makeatletter
 \theoremstyle{plain}    
 \newtheorem{thm}{Theorem}[section]
 \numberwithin{equation}{section} 
 \numberwithin{figure}{section} 
 \theoremstyle{plain}
 \theoremstyle{plain}    
 \newtheorem{prop}[thm]{Proposition} 
 \theoremstyle{plain}    
 \newtheorem{lem}[thm]{Lemma} 
 \theoremstyle{plain}    
 \newtheorem{cor}[thm]{Corollary} 
 \theoremstyle{remark}
 \newtheorem{rem}[thm]{Remark}
 \theoremstyle{remark}    
 \newtheorem{acknowledgement}[thm]{Acknowledgement} 

\def\makebbb#1{
    \expandafter\gdef\csname#1\endcsname{
        \ensuremath{\Bbb{#1}}}
}
\makebbb{R}
\makebbb{N}
\makebbb{Z}
\makebbb{C}
\makebbb{G}
\makebbb{E}
\makebbb{H}
\makebbb{P}
\makebbb{B}
\makebbb{K}

\usepackage{babel}
\makeatother
\begin{document}

\title{Super Toeplitz Operators on Line Bundles}

\author{Robert Berman}

\email{robertb@math.chalmers.se}


\begin{abstract}
Let $L^{k}$ be a high power of a hermitian holomorphic line bundle
over a complex manifold $X.$ Given a differential form $f$ on $X,$
we define a super Toeplitz operator $T_{f}$ acting on the space of
harmonic $(0,q)$-forms with values in $L^{k},$ with symbol $f.$
The asymptotic distribution of its eigenvalues, when $k$ tends to
infinity, is obtained in terms of the symbol of the operator and the
curvature of the line bundle $L$, given certain conditions on the
curvature. For example, already when $q=0,$ i.e. the case of holomorphic
sections, this generalizes a result of Boutet de Monvel and Guillemin
to semi-positive line bundles. The asympotics are obtained from the
asymptotics of the Bergman kernels of the corresponding harmonic spaces.
Applications to sampling are also given. 
\end{abstract}
\maketitle

\section{Introduction}

Let $(X,\omega)$ be an $n-$dimensional compact hermitian manifold
and let $L$ be a hermitian holomorphic line bundle over $X.$ The
fiber metric on $L$ will be denoted by $\phi.$ It can be thought
of as a collection of local functions: let $s$ be a local holomorphic
trivializing section of $L,$ then locally, $\left|s(z)\right|_{\phi}^{2}=e^{-\phi(z)}$
and the canonical curvature two-form of $L$ is $\partial\overline{\partial}\phi.$
Denote by $X(q)$ be the subset of $X$ where the curvature two-form
of $L$ is non-degenerate and has exactly $q$ negative eigenvalues.
The notation $\eta_{p}:=\eta^{p}/p!$ will be used in the sequel,
so that the volume form on $X$ may be written as $\omega_{n}.$

The spaces $H^{0}(X,L^{k}),$ consisting of global holomorphic sections
with values in high powers of $L,$ appear naturally in complex and
algebraic geometry, as well as in mathematical physics. In many applicaions
the line bundle $L$ is positive i.e. its curvature two-form is positive
and the asymptotic properties of the sequence of Hilbert spaces $H^{0}(X,L^{k})$
have been studied thorougly in this case. For example the asymptotic
behaviour of the corresponding Bergman kernels is known and can be
used to study asymptotic properties of Toeplitz operators acting on
$H^{0}(X,L^{k})$ as well as asymptotic conditions on the density
of the distribution of sampling points on the manifold $X$ (see \cite{bern}
for a recent survey from this point of view). The aim of the present
article is to extend these results in two directions: to line bundles
with weaker curvature properties than positivity, such as semi-positivity
(part $1$) and to harmonic $(0,q)-$ forms with values in $L^{k}$
(part $2$). To emphasize the analogy between holomorphic sections
and harmonic forms, some rudiments of the theory of super manifolds
is recalled. The super formalism also offers a compact notation. 

In part $1,$ following \cite{li} and \cite{bern}, everything is
reduced to knowing the leading asymptotic behaviour of the Bergman
kernel $K(x,y).$ The asymptotics, in turn, are obtained using a new
and comparatively elementary approach based on the method used in
\cite{berm} to prove local holomorphic Morse inequalities. The main
application is a generalization of a theorem of Boutet de Monvel and
Guillemin that expresses the asymptotic distribution of the eigenvalues
of a Toeplitz operator in terms of the symbol of the operator \cite{bo-gu},
\cite{gu}. When $L$ is positive the associated dual disc bundle
over $X$ is strictly pseudoconvex. One can then profit from the knowledge
of the Bergman kernel on a strictly pseudoconvex manifold \cite{bo-sj}.
However, when $L$ is only semi-positive one would have to use the
corresponding result on a weakly pseudoconvex manifold, which is not
available. In fact, a recent counter example of Donnelly \cite{do}
to a conjecture due to Siu, shows that the tangential Cauchy-Riemann
operator on the boundary of the dual disc bundle does not have closed
range. This property is essential to the previous approaches to the
asymptotics of the Bergman kernel. 

In part $2$ the approach in part $1$ is extended to study the Bergman
kernel of the space of harmonic $(0,q)-$ forms with values in $L^{k},$
considered as a bundle valued form on $X\times X.$ The main application
is a generalization of the theorem of Boutet de Monvel and Guillemin
to Toeplitz operators, whose symbol is a differential form on $X.$
These operators are called super Toeplitz operators and they are closely
related to the operators introduced in \cite{bor} in the context
of Berezin-Toeplitz quantization of symplectic super manifolds.

It should be added that part one is just a special case of part two
(when $q$ is equal to zero), except for the applications to sampling.
However, it has been included to motivate the more general discussion
given in the second part.

\part{Holomorphic sections}

Let $(\psi_{i})$ be an orthormal base for $H^{0}(X,L).$ Denote by
$\pi_{1}$ and $\pi_{2}$ the projections on the factors of $X\times X.$
The \emph{Bergman kernel} of the Hilbert space $H^{0}(X,L)$ is defined
by \[
K(x,y)=\sum_{i}\psi_{i}(x)\otimes\overline{\psi_{i}(y)}.\]
 Hence, $K(x,y)$ is a section of the pulled back line bundle $\pi_{1}^{*}(L)\otimes\overline{\pi_{2}^{*}(L)}$
over $X\times X.$ For a fixed point $y$ we identify $K_{y}(x):=K(x,y)$
with a section of the hermitian line bundle $L\otimes L_{y},$ where
$L_{y}$ denotes the line bundle over $X,$ whose constant fiber is
the fiber of $L$ over $y,$ with the induced metric. The definition
of $K$ is made so that $K$ satifisfies the following reproducing
property \begin{equation}
\alpha(y)=(\alpha,K_{y})\label{(I)repr property}\end{equation}
\footnote{We are abusing notation here: the scalar product $(\cdot,\cdot)$
on $H^{0}(X,L)$ determines a pairing of $K_{y}$ with any element
of $H^{0}(X,L),$ yielding an element of $L_{y}.$ %
}for any element $\alpha$ of $H^{0}(X,L),$ which also shows that
$K$ is well-defined. In other words $K$ represents the orthogonal
projection onto $H^{0}(X,L)$ in $L^{2}(X,L).$ The restriction of
$K$ to the diagonal is a section of $L\otimes\overline{L}$ and we
let $B(x)=\left|K(x,x)\right|$ be its pointwise norm: \[
B(x)=\sum_{i}\left|\psi_{i}(x)\right|^{2}.\]
We will refer to $B(x)$ as the \emph{Bergman function} of $H^{0}(X,L).$
It has the following extremal property:\begin{equation}
B(x)=\sup\left|\alpha(x)\right|^{2},\label{(I)extremal prop of B}\end{equation}
 where the supremum is taken over all normalized elements $\alpha$
of $H^{0}(X,L).$ An element realizing the extremum, is called an
\emph{extremal at the point $x$} and is determined up to a complex
constant of unit norm. In order to estimate $K(x,y)$ we will have
great use for a more general identity. It is just a reformulation
of the fact that, by the reproducing property \ref{(I)repr property},
$K_{x}/\sqrt{B(x)}$ may be identified with an extremal at the point
$x.$ 

\begin{prop}
\label{(I) prop: K is extremal}Let $\alpha$ be an extremal at the
point $x.$ Then\[
\left|K(x,y)\right|^{2}=\left|\alpha(y)\right|^{2}B(x)\]

\end{prop}
\begin{proof}
First fix the point $x$ and take a local holomorphic trivialization
of $L$ around $x.$ Then we may identify $K(x,y)$ with an element
$K_{x}$ of $H^{0}(X,L).$ Now we may assume that $\left\Vert K_{x}\right\Vert \neq0$
- it will be clear that otherwise the statement is trivially true.
By the reproducing property \ref{(I)repr property} of $K$ the normalized
element $K_{x}/\left\Vert K_{x}\right\Vert $ is an extremal of $B$
at $x.$ Furthermore, the reproducing property \ref{(I)repr property}
also shows that the squared norm of $K_{x}$ is given by $(K_{x},K_{x})=K_{x}(x)=K(x,x).$
Hence \[
\left|\alpha(y)\right|^{2}=\left|K_{x}(y)\right|^{2}/K(x,x),\]
 for any other extremal $\alpha$ of $B$ at $x.$ Since $\left|K(x,y)\right|^{2}=\left|K_{x}(y)\right|^{2}e^{-\phi(x)}$
and since by definition $K(x,x)e^{-\phi(x)}=B(x)$ this proves the
proposition. 
\end{proof}
Next, we will define certain operators on $H^{0}(X,L).$ Given a complex-valued
bounded mesurable function $f$ on $X$ we define $T_{f},$ the so
called \emph{Toeplitz operator with symbol $f,$} by \[
T_{f}:=P\circ f\cdot,\]
 where $f\cdot$ denotes the usual multiplication operator on $L^{2}(X,L)$
and $P$ is the orthogonal projection onto $H^{0}(X,L).$ Equivalently:
\begin{equation}
(T_{f}\alpha,\beta)=(f\alpha,\beta),\label{(I) Toplitz using scalar pr}\end{equation}
 for all elements $\alpha$ and $\beta$ of $H^{0}(X;L).$ Note that
the operator $T_{f}$ is hermitian if $f$ is real-valued.

When studying asymptotic properties of $L^{k},$ all objects introduced
above will be defined with respect to the line bundle $L^{k}.$

\section{\label{(I) Section: Asym for B and T}Asymptotic results for Bergman
kernels and Toeplitz operators.}

Let us first see how to prove the following upper bound on the Bergman
kernel function $B(x):$ \begin{equation}
B^{k}(x)\leq k^{n}\frac{1}{\pi^{n}}1_{X(0)}(x)\left|det_{\omega}(\frac{i}{2}\partial\overline{\partial}\phi)_{x}\right|+o(k^{n}),\label{(I) upper bound on B}\end{equation}
 where we have identified the two-form $i\partial\overline{\partial}\phi$
with an endomorphism, using the metric $\omega$, so that its determinant
is well-defined. Integrating this over all of $X$ gives an upper
bound on the dimension of the space of holomorphic sections:\begin{equation}
\dim_{\C}H_{\overline{\partial}}^{0}(X,L^{k})\leq k^{n}\frac{(-1)^{q}}{\pi^{n}}\frac{1}{n!}\int_{X(0)}(\frac{i}{2}\partial\overline{\partial}\phi)^{n}+o(k^{n}),\label{(I) upper bound on H0}\end{equation}
 which is precisely Demailly's holomorphic Morse-inequalities for
$(0,q)-$forms when $q=0.$ In \cite{berm} the inequality \ref{(I) upper bound on B}
and it's generalization to harmonic $(0,q)-$ forms were called \emph{local}
holomorphic Morse-inequalties. As we will see \ref{(I) upper bound on B}
follows from the submean property of holomorhic functions and a simple
localization argument. Fix a point $x$ in $X$ and choose complex
coordinates $z$ and a holomorphic trivialization $s$ of $L$ around
$x$, such the metric $\omega$ is Euclidean with respect to $z$
at $0$ and the fiber metric $\phi(z)=\phi_{0}(z)+O(\left|z\right|^{3}),$
where $\phi_{0}(z)=\sum_{i=1}^{n}\lambda_{i}\left|z_{i}\right|^{2}$
and $\lambda_{i}$ are the eigen values of the curvature two-form
$\partial\overline{\partial}\phi,$ with respect to the base metric
$\omega,$ at the point $x.$ According to the extremal property \ref{(I)extremal prop of B}
of $B(x)$ we have to estimate the pointwise norm of section $\alpha$
at $x$ in terms of the global $L^{2}$ norm. Let $B_{R_{k}}$ be
balls centered at $x$ of radius $R_{k}\rightarrow0.$ By first restricting
the global norm to the ball $B_{R_{k}}$ and than making the change
of variables $z=\frac{w}{\sqrt{k}}$ in the integral we get \[
\frac{\left|\alpha_{k}(x)\right|^{2}}{\left\Vert \alpha_{k}\right\Vert _{X}^{2}}\leq\frac{\left|\alpha_{k}(x)\right|^{2}}{\left\Vert \alpha_{k}\right\Vert _{R_{k}}^{2}}=\rho_{k}k^{n}\left|f_{k}(0)\right|^{2}/\int_{B_{\sqrt{k}R_{k}}}\left|f_{k}(w)\right|^{2}e^{-\phi_{0}(w)},\]
 where the holomorphic functions $f_{k}$ represent $\alpha_{k}$
in the local frame. The factor $\rho_{k}$ comes from the base manifold
metric and the terms of order $O(\left|z\right|^{3})$ in the fiber
metric on $L$. If we now choose e.g. $R_{k}=\frac{lnk}{\sqrt{k}}$
then the factor $\rho_{k}\rightarrow1$ and the scaled radii $\sqrt{k}R_{k}\rightarrow\infty$
so that the integration in the variable $w$ is over all of $\C^{n}$
in the limit. Furthermore, since $\left|f_{k}\right|^{2}$ is plurisubharmonic
the quotient in the right hand side can be estimated by the inverse
of the Gaussian \[
1/\int_{B_{\sqrt{k}R_{k}}}e^{-\phi_{0}(w)}\]
 wich tends to $(1/\pi)^{n}\lambda_{1}\lambda_{2}\cdot\cdot\cdot\lambda_{n}$
if all eigenvalues are positive and is equal to zero in the limit
otherwise. This proves the upper bound on the Bergman kernel function
\ref{(I) upper bound on B}. In fact, what we have proved is the stronger
statement that for any sequence $(\alpha_{k}),$ where $\alpha_{k}$
is in $H^{0}(X,L^{k}),$\begin{equation}
\limsup_{k}k^{-n}\left|\alpha_{k}(x)\right|^{2}/\left\Vert \alpha_{k}\right\Vert _{B_{R_{k}}}^{2}\leq\frac{1}{\pi^{n}}1_{X(0)}(x)\left|det_{\omega}(\frac{i}{2}\partial\overline{\partial}\phi)_{x}\right|\label{(I) local morse q is 0}\end{equation}
 which will be important in the proof of theorem \ref{(I) thm: asym for K}.

It is well-known that \ref{(I) upper bound on B} is actually an asymptotic
\emph{equality} when $L$ is positive on all of $X.$ One can e.g.
use Hörmander's celebrated $L^{2}$ estimates to obtain the equality
\cite{ti},\cite{li}, \cite{bern} (a complete asymptotic expansion
is obtained in \cite{ze}, using micro local analysis based on \cite{bo-sj}).
But these methods break down if the curvature of $L$ is only semi-positive.
On the other hand Demailly proved, using his holomorphic Morse inequalitues,
that \ref{(I) upper bound on B} is an asymptotic equality under the
more general condition that $X(1)$ is empty. Combining Demailly's
result with the upper bound \ref{(I) upper bound on B} we obtain
the following theorem:

\begin{thm}
\label{(I) Thm: asym for B}Suppose that $X(1)$ is empty. Then\[
k^{-n}B^{k}(x)\rightarrow\frac{1}{\pi^{n}}1_{X(0)}(x)\left|det_{\omega}(\frac{i}{2}\partial\overline{\partial}\phi)_{x}\right|\]
in $L^{1}(X,\omega_{n}).$ In particular, the measure $B^{k}\omega_{n}/k^{n}$
converges to $\pi^{-n}1_{X(0)}(\partial\overline{\partial}\phi)_{n}$
in the weak{*}-topology. 
\end{thm}
\begin{proof}
The upper bound \ref{(I) upper bound on B} says that

\[
\limsup_{k}\begin{array}{lr}
k^{-n}B^{k}(x)\leq\frac{1}{\pi^{n}}1_{X(0)}(x)\left|det_{\omega}(\frac{i}{2}\partial\overline{\partial}\phi)_{x}\right|\end{array},\]
 for any line bundle $L.$ Moreover, if the curvature of the line
bundle $L$ is such that $X(1)$ is empty, then \[
\lim_{k}k^{-n}\int_{X}B^{k}(x)\omega_{n}=\frac{1}{\pi^{n}}\int1_{X(0)}(x)\left|det_{\omega}(\frac{i}{2}\partial\overline{\partial}\phi)_{x}\right|\omega_{n}.\]
To see this, note that the left hand side is the dimension of the
space $H^{0}(X,L^{k}).$ In this form the statement was first shown
by Demailly in \cite{d1}. See also proposition \ref{(II) Prop: equality for Hq}
in the present paper. Finally the theorem follows from the following
simple lemma:
\end{proof}
\begin{lem}
\label{(I) lemma: int theory}Assume that $(X,\mu)$ is a finite measure
space and that $f$ and $f_{k}$ are bounded functions, where the
sequence $f_{k}$ is uniformly bounded. If \[
(\textrm{i})\lim_{k}\int_{X}f_{k}d\mu=\int_{_{X}}fd\mu,\,\,\,\,\textrm{and}\,\,\,\,(\textrm{ii})\,\limsup f_{k}\leq f.\]
Then the sequence $f_{k}$ converges to $f$ in $L^{1}(X,\mu).$ 
\end{lem}
\begin{proof}
By the assumption (i)\[
\limsup_{k}\int_{X}\left|f_{k}-f\right|d\mu=2\limsup_{k}\int_{X}\chi_{k}(f_{k}-f)d\mu,\]
 where $\chi_{k}$ is the characteristic funtion of the set where
$f_{k}-f$ is non-negative. The right hand side can be estimated by
Fatou's lemma, which is equivalent to the inequality \[
\limsup_{k}\int_{X}g_{k}d\mu\leq\int_{X}\limsup_{k}g_{k}d\mu,\]
 if the sequence $g_{k}$ is dominated by an $L^{1}-$function. Taking
$g_{k}=\chi_{k}(f_{k}-f)$ and using the assumption $(ii),$ finishes
the proof of the lemma.
\end{proof}
The weak convergence of the previous theorem can be reformulated in
terms of Toeplitz operators:

\begin{cor}
\label{(I) cor: Tr T}Suppose that $X(1)$ is empty. Then for any
bounded function $f$ on $X$
\end{cor}
\[
\lim_{k}\textrm{k}^{-n}\textrm{Tr}T_{f}=(2\pi)^{-n}\int_{X(0)}f(i\partial\overline{\partial}\phi)_{n}\]

\begin{proof}
From the definition \ref{(I) Toplitz using scalar pr} of a Toeplitz
operator $\textrm{Tr}T_{f}=\sum_{i}(f\Psi_{i},\Psi_{i}),$ which is
equal to $\int_{X}fB_{k}(x)\omega_{n}.$ The corollary now follows
from the $L^{1}$ convergence in the previous theorem.
\end{proof}
We have seen how to obtain the leading asymptotics of $B(x),$ the
norm of the Bergman kernel $K$ on the diagonal.The main point of
the present paper is that the argument presented actually also shows
that $\left|K(x,y)\right|^{2}$ tends to zero off the diagonal to
the leading order. The idea is that combining the upper bound \ref{(I) local morse q is 0}
with the asymptotics for $B(x),$ one sees that any sequence of extremals
$\alpha_{k}$ at a given point $x,$ becomes localized around $x$
in the large $k$ limit. Since $K_{x}(y)$ is essentually equal to
$\alpha_{k}(y),$ this will show that $K(x,y)$ localizes to the diagonal
when $k$ tends to infinity. A similar argument has been used by Bouche
\cite{bouche} to construct holomorphic peak sections when $L$ is
positive.

\begin{thm}
\label{(I) thm: asym for K}Suppose that $X(1)$ is empty. Denote
by $\Delta$ the current of integration on the diagonal in $X\times X.$
Then \[
\begin{array}{lr}
\lim_{k}k^{-n}\left|K^{k}(x,y)\right|^{2}\omega_{n}(x)\wedge\omega_{n}(y)=(2\pi)^{-n}1_{X(0)}\Delta\wedge(i\partial\overline{\partial}\phi)_{n}\end{array},\]
as measures on $X\times X$, in the weak {*}-topology. 
\end{thm}
\begin{proof}
First note that the mass of the measures $\mu_{k}:=\left|K^{k}(x,y)\right|^{2}/k^{n}\omega_{n}(x)\wedge\omega_{n}(y)$
are uniformly bounded in $k$: first integrating over $y$ and using
the reproducing property \ref{(I)repr property} gives \[
\mu_{k}(X\times X)=k^{-n}\int_{X_{x}}\left|(K_{x}^{k},K_{x}^{k})\right|\omega_{n}(x)=k^{-n}\int_{X_{x}}B_{k}(x)\omega_{n}=k^{-n}\dim H^{0}(X,L^{k}),\]
 which clearly is bounded by \ref{(I) upper bound on H0}. Moreover,
the mass of $X(0)^{c}\times X$ tends to zero: \[
\mu_{k}(X(0)\times X)=k^{-n}\int_{X(0)}\left|(K_{x}^{k},K_{x}^{k})\right|=k^{-n}\int_{X(0)}B_{k}(x)\omega_{n}=0,\]
 where we have used the bound \ref{(I) upper bound on B}.

Hence it is enough to prove the convergence on $X(0)\times X$. Moreover
since the mass of the measures $\mu_{k}$ is bounded, it is enough
to show that any subsequnce of $\mu_{k}$ has another subsequence
that converges to $\pi^{-n}1_{X(0)}\Delta\wedge(\partial\overline{\partial}\phi)_{n}$,
in the weak{*}-topology. To simplify the notation the first subsequence
will be indexed by $k$ in the following.

According to theorem \ref{(I) Thm: asym for B} and standard integration
theory there is a subsequence of $k^{-n}B_{k}$ that converges to
$\pi{}^{-n}\left|det_{\omega}(\frac{i}{2}\partial\overline{\partial}\phi)_{x}\right|$
almost everywhere on $X(0).$ Fix a point $x$ in $X(0)$ where $k^{-n}B_{k}(x)$
converges. Take a sequence of sections $\alpha_{k},$ where $\alpha_{k}$
is a normalized extremal at the point $x.$ Then according to proposition
\ref{(I) prop: K is extremal}: \begin{equation}
\left|K_{x}^{k}(y)\right|^{2}=\left|\alpha_{k}(y)\right|^{2}B^{k}(x)\label{(I) K is extremal}\end{equation}
We will now show that there is a subsequence of $\alpha_{k}$ such
that \begin{equation}
\lim_{k}\left\Vert \alpha_{k}\right\Vert _{B_{R_{k}}(x)}^{2}=1\label{(I) claim: peak sections}\end{equation}
 where $B_{R_{k}}(x)$ is a ball centered in $x$ of radius $k$ of
radius $R_{k}:=\textrm{l}nk/\sqrt{k}$ with respect to the {}``normal''
coordinates around $x$ used in section \ref{(I) Section: Asym for B and T}.
Since the global norm of $\alpha_{k}$ is equal to one and the radii
$R_{k}$ tend to zero, it follows that the function $\left|\alpha_{k}(y)\right|^{2}$
on $X$ converges to the Dirac measure at $x$ in the weak{*}-topology.
From this convergence we will be able to deduce the statement of the
theorem. To prove the claim \ref{(I) claim: peak sections} first
observe that there is a subsequence of $\alpha_{k}$ such that at
the point $x:$ \[
\pi{}^{-n}\left|det_{\omega}(\frac{i}{2}\partial\overline{\partial}\phi)_{x}\right|=\lim_{j}k_{j}^{-n}\left|\alpha_{k_{j}}(x)\right|^{2}\]
Indeed, since $\alpha_{k}$ is normalized \ref{(I)extremal prop of B}
says that the right hand side is equal to $k^{-n}$ times $B_{k}(x),$
the Bergman function at the point $x,$ which in turn tends to the
left hand side according to theorem \ref{(I) Thm: asym for B} and
by the assumption on the point $x.$ 

Furthermore, since $\alpha_{k}$ is normalized the restricted norm
$\left\Vert \alpha_{k}\right\Vert _{B_{R_{k}}(x)}^{2}$ is less than
one. Hence the right hand side is trivially estimated by \[
\limsup_{j}k_{j}^{-n}\left|\alpha_{k_{j}}(x)\right|^{2}/\left\Vert \alpha_{k_{j}}\right\Vert _{B_{R_{k_{j}}}(x)}^{2}.\]
 According to \ref{(I) local morse q is 0} this in turn may be estimated
by $\left|det_{\omega}(\frac{i}{2}\partial\overline{\partial}\phi)_{x}\right|.$
All in all this shows that \[
\pi{}^{-n}\left|det_{\omega}(\frac{i}{2}\partial\overline{\partial}\phi)_{x}\right|=\pi{}^{-n}\lim_{j}\left|det_{\omega}(\frac{i}{2}\partial\overline{\partial}\phi)_{x}\right|/\left\Vert \alpha_{k_{j}}\right\Vert _{B_{R_{k_{j}}}(x)}^{2}\]
 Since the left hand side is non-zero on $X(0)$ this proves the claim
\ref{(I) claim: peak sections}. Now to prove the theorem we take
a test function $f(x,y)$ and consider the integral \[
k_{j}^{-n}\int_{X(0)\times X}f(x,y)\left|K^{k_{j}}(x,y)\right|^{2}\omega_{n}(x)\wedge\omega_{n}(y).\]
We may restrict the integration over the first factor in $X(0)\times X$
to the set of all $x$ satisfying the assumption used above, since
the the complement is of mesure zero. Using the identity \ref{(I) K is extremal}
the integral over $y$ (for a fixed point $x)$ equals \[
k_{j}^{-n}B^{k_{j}}(x)\int_{X}f(x,y)\left|\alpha_{k_{j}}(y)\right|^{2}\omega_{n}(y).\]
Since by \ref{(I) claim: peak sections} the function $\left|\alpha_{k}(y)\right|^{2}$
converges to the Dirac measure at $x$ in the weak{*}-topology, the
integral is equal to $f(x,x)$ in the limit and the first factor is
equal to $\pi{}^{-n}\left|det_{\omega}(\frac{i}{2}\partial\overline{\partial}\phi)_{x}\right|$
in the limit (by theorem \ref{(I) Thm: asym for B}). Hence the previous
integral is equal to \[
(2\pi)^{-n}\int_{X(0)}f(x,x)(i\partial\overline{\partial}\phi)_{n}\]
 in the limit, which finishes the proof of the theorem.
\end{proof}
As before the convergence may be formulated in terms of Toeplitz operators
(a more general statement will be proved in part II (corollary \ref{(II) cor: tr TT})).

\begin{cor}
\label{(I) cor: Tr TT}Suppose that $X(1)$ is empty. Then \[
\begin{array}{lr}
\lim_{k}\textrm{k}^{-n}\textrm{Tr}\, T_{fg}=\lim_{k}k^{-n}\textrm{Tr}\, T_{f}T_{g}\end{array}.\]

\end{cor}
We will now use the results on the asymptotics of the Bergman kernel
$K$ to express asymptotic spectral properties of Toeplitz operators
in terms of their symbol. Denote by \[
N(T_{f}>\gamma)\]
 the number of eigenvalues of $T_{f}$ that are greater than the number
$\gamma$ (counted with multiplicity). Furthermore, $N(T_{f}<\gamma)$
is defined similarly. 

\begin{thm}
\label{(I) thm: N is}Suppose that $X(1)$ is empty and that $f$
is a real-valued bounded function. Then for all $\gamma$ except possibly
countably many the following holds:

\begin{equation}
\lim_{k}k^{-n}N(T_{f}>\gamma)=(2\pi)^{-n}\int_{\left\{ f>\gamma\right\} \bigcap X(0)}(i\partial\overline{\partial}\phi)_{n}\label{(I) statement of thm N is}\end{equation}
 and similarly for $N(T_{f}<\gamma).$ 
\end{thm}
\begin{proof}
Given the asymptotic behaviour of $K(x,y)$ in theorem \ref{(I) thm: asym for K},
the proof can be adapted word by word from \cite{li},\cite{bern}.
But for completeness we give a proof here, that slightly simplifies
the proof in \cite{li}. We first prove the statement when $f$ is
the characterstic function for a given set $\Omega:$ $f=1_{\Omega}.$
Let us denote by $T_{\Omega}$ the corresponding operator. We may
assume that $1>\gamma>0.$ By corollary \ref{(I) cor: Tr T} the right
hand side of \ref{(I) statement of thm N is} is then equal to the
limit of $k^{-n}\textrm{Tr T}_{\Omega}.$ Moreover, by corollary \ref{(I) cor: Tr TT}
this limit in turn is equal to the limit of $k^{-n}\textrm{TrT}_{\Omega}^{2}.$
We will now see that this can happen only if $\lim_{k}k^{-n}N(T_{\Omega}>\gamma)=\lim_{k}k^{-n}\textrm{Tr T}_{\Omega},$
which proves the statement with this special choice of $f.$ Indeed,
since if we denote by $\tau_{j}$ the eigen values of $T_{\Omega},$
\[
\lim_{k}k^{-n}\sum\tau_{j}(1-\tau_{j})=0,\]
 it follows from estimating $(1-\tau_{j})$ from below that \[
\lim_{k}k^{-n}\sum_{\tau_{j}\leq\gamma}\tau_{j}=0.\]
 Hence\[
\lim_{k}k^{-n}\sum_{\tau_{j}>\gamma}\tau_{j}=\lim_{k}k^{-n}\textrm{Tr T}_{\Omega}=\int_{\Omega\bigcap X(0)}(\partial\overline{\partial}\phi)_{n}\]
 Now it is not hard to deduce that $\lim_{k}k^{-n}N(T_{\Omega}>\gamma)=\int_{\Omega\bigcap X(0)}(\partial\overline{\partial}\phi)_{n}.$ 

By comparing an arbitrary function $f$ with a characteristic function
we will now finish the proof of the theorem. Let us first prove the
lower bound \begin{equation}
\liminf_{k}k^{-n}N(T_{f}>\gamma)\geq(2\pi)^{-n}\int_{\left\{ f>\gamma\right\} \bigcap X(0)}(i\partial\overline{\partial}\phi)_{n}.\label{(I) pf of N is: lower bd}\end{equation}
First note that we may assume that $f$ is non-negative by adding
an appropriate constant to $f.$ By the max-min principle applied
to the operator $T_{f}$ and by \ref{(I) Toplitz using scalar pr}
\[
N(T_{f}>\gamma)=\max\left\{ \dim V:\,(f\alpha,\alpha)>\gamma(\alpha,\alpha)\,\forall\alpha\in V\right\} ,\]
 where $V$ is a linear subspace of $H^{0}.$ Hence we have to find
a sequence of subspaces $V_{k}$ with \begin{equation}
\dim V_{k}=k^{n}(2\pi)^{-n}\int_{\left\{ f>\gamma\right\} \bigcap X(0)}(i\partial\overline{\partial}\phi)_{n}+o(k^{n}),\label{(I) pf of N is: dim}\end{equation}
 such that for any normalized $\alpha$ in $V_{k}$ $(f\alpha,\alpha)>\gamma.$
To this end denote by $\Omega$ the set where $f>\gamma$ and denote
by $\chi$ the corresponding characteristic function. Since we have
already proved the theorem for characteristic functions, there is,
for any given small positive $\varepsilon,$ a sequence of subspaces
$V_{k}$ with the correct dimension \ref{(I) pf of N is: dim} such
that \[
(\chi\alpha,\alpha)>1-\varepsilon\]
 for all $\alpha$ in $V_{k},$ Since by definition $f>\gamma\chi$
it follows that \[
(f\alpha,\alpha)>\gamma(1-\varepsilon),\]
 for all $\alpha$ in $V_{k}.$ By symmetry this means that \[
\liminf_{k}k^{-n}N(T_{f}>\gamma)\geq(2\pi)^{-n}\int_{\left\{ f>\gamma+\varepsilon\right\} \bigcap X(0)}(i\partial\overline{\partial}\phi)_{n}.\]
 By letting $\epsilon$ tend to zero in the right hand side we obtain
the desired lower bound \ref{(I) pf of N is: lower bd}. If we now
apply this result to the function $-f$ we obtain the following equivalent
result: \begin{equation}
\liminf_{k}k^{-n}N(T_{f}<\delta)\geq(2\pi)^{-n}\int_{\left\{ f<\delta\right\} \bigcap X(0)}(i\partial\overline{\partial}\phi)_{n}.\label{(I) lower bd two}\end{equation}
Now note that for all except countably many numbers $\gamma,$ the
set $\{ f=\gamma\}$ is of mesure zero with respect to the mesure
$(i\partial\overline{\partial}\phi)_{n}$ on $X.$ Indeed, the function
$g(\gamma)=\int_{\left\{ f\leq\gamma\right\} \bigcap X(0)}(i\partial\overline{\partial}\phi)_{n}$
on the real line is increasing and it is well-known that an increasing
function is continous accept on a countable set. This forces the mesure
of the set $\{ f=\gamma\}$ to be zero for all $\gamma$ accept those
in the countable set. Finally, since the total number of eigen values
for any operator $T_{f}$ is equal to the dimension of $H^{0}(X,L^{k}),$
we get, when $k$ tends to infinity, that the sum\[
\lim_{k}k^{-n}(N(T_{f}>\gamma)+\lim_{k}k^{-n}N(T_{f}\leq\gamma)\]
is equal to the asymptotic dimension\[
(2\pi)^{-n}(\int_{\left\{ f>\gamma\right\} \bigcap X(0)}(i\partial\overline{\partial}\phi)_{n}+(2\pi)^{-n}\int_{\left\{ f<\gamma\right\} \bigcap X(0)}(i\partial\overline{\partial}\phi)_{n})\]
 for all $\gamma$ such that the mesure of the set $\{ f=\gamma\}$
is zero. Combing this with the lower bounds \ref{(I) pf of N is: lower bd}
and \ref{(I) lower bd two}, we see that we must have equality in
\ref{(I) pf of N is: lower bd}, which proves the theorem. 
\end{proof}
The main application of the previous theorem is to show that there
is a large supply of holomorphic sections concentrated on any given
set $\Omega$ in $X(0),$ in the following sence: \[
\left\Vert \alpha\right\Vert _{\Omega}^{2}\geq(1-\varepsilon)\left\Vert \alpha\right\Vert _{X}^{2}\]
for any given positive (small) $\epsilon.$ To see this, denote by
$\chi$ the characteristic function of the set $\Omega,$ and note
that if $\alpha$ is a linear combination of eigen sections of the
Toeplitz operator $T_{\chi}$ then $\alpha$ will be concentrated
on $\Omega,$ as long as the eigenvalues are bounded from below by
$(1-\varepsilon).$ The number of such sections $\alpha$ is precisely
the spectral counting function $N(T_{\chi}\geq1-\varepsilon).$ Hence,
the previous theorem shows that there is a subspace of dimension \[
k^{n}(2\pi)^{-n}\int_{\Omega}(i\partial\overline{\partial}\phi)_{n}+o(k^{n}),\]
consisting of concentrated sections, a result that will be useful
when studying sampling sequences in the next section. 

The following equivalent formulation of theorem \ref{(I) thm: N is}
can be obtained by standard methods in spectral theory. It generalizes
a theorem of Boutet de Monvel and Guillemin \cite{bo-gu}, \cite{gu},
valid when $L$ is positive, to the case when $X(1)$ is empty.

\begin{thm}
Suppose that $X(1)$ is empty. Let $(\tau_{i})$ be the eigen values
of $T_{f}$ and denote by $d\xi_{k}$ the spectral measure of $T_{f}$
divided by $k^{n},$ i.e. \[
d\xi_{k}:=k^{-n}\sum_{i}\delta_{\tau_{i}},\]
 where $\delta_{\tau_{i}}$ is the Dirac measure centered at $\tau_{i}.$
Then $d\xi_{k}$ tends, in the weak{*}-topology, to the push forward
of the measure $1_{X(0)}(2\pi)^{-n}(i\partial\overline{\partial}\phi)_{n}$
under the map $f,$ i.e. \[
\lim_{i}\sum_{i}a(\tau_{i})=(2\pi)^{-n}\int_{X(0)}a(f(x)(i\partial\overline{\partial}\phi)_{n}\]
 for any mesurable function $a$ on the real line.
\end{thm}

\section{Sampling}

Let $D_{k}$ be a finite set of points in $X.$ We say that the sequence
of sets $D_{k}$ is \emph{sampling} for the sequence of Hilbert spaces
$H^{0}(X,L^{k})$ if there exists a uniform constant $A$ such that
\[
A^{-1}k^{-n}\sum_{D_{k}}\left|\alpha(x)\right|^{2}\leq\left\Vert \alpha\right\Vert ^{2}\leq Ak^{-n}\sum_{D_{k}}\left|\alpha(x)\right|^{2},\]
 for any element $\alpha$ in $H^{0}(X,L^{k}).$ The points in $D_{k}$
will assumed to be \emph{separated} in the following sense: the distance
between any two points in $D_{k}$ is bounded from below by a uniform
constant times $k^{-n}.$ Consider the measures $d\nu_{k}$ on $X$
corresponding to the sets $D_{k}:$ \[
d\nu_{k}:=k^{-n}\sum_{D_{k}}\delta_{x}.\]
 Because of the separability assumption their mass is uniformly bounded
in $k.$ Hence any subsequence has a subsequence that is weak {*}-convergent.
Denote by $d\nu$ such a limit measure. It is natural to ask how dense
the sampling points should be, for large $k,$ in order to be sampling.
i.e. we ask for asymptotic density conditions on the measure $d\nu.$
The model case is sampling on lattices for the Fock space, i.e. $X$
is taken to be $\C^{n}$ with its standard Euclidean metric form $\omega$and
$L$ is the line bundle with constant positive curvature $-2i\omega.$
If the sequence $D_{k}$ is a sequence of lattices generated over
$\Z$ by $k^{-1/2}(a_{1},...,a_{2n}),$ where the $a_{i}$ are positive
numbers, then a necessary condition for this sequence to be sampling
is that $a_{1}\cdot\cdot\cdot a_{2n}\leq\pi^{n}.$ In (\cite{li},
\cite{bern}) this necessary condition was generalized to any positive
line bundle. Namely, the limit measure has to satisfy \[
d\nu\geq(2\pi)^{-n}(i\partial\overline{\partial}\phi)_{n}.\]

The next theorem shows that in order to sample $H^{0}(X,L^{k})$ when
$X(1)$ is empty, the sampling points have to satisfy the same necessary
density conditions in $X(0)$ (the part of $X$ where $L$ is positive)
as in the case when the curvature is positive everywhere on $X.$

\begin{thm}
Assume that the sequence of sets $D_{k}$ is sampling for the sequence
of Hilbert spaces $H^{0}(X,L^{k}).$ If $X(1)$ is empty then the
following necessary condition holds: \[
d\nu\geq(2\pi)^{-n}(i\partial\overline{\partial}\phi)_{n}\]
 on $X(0).$ In other words \begin{equation}
\liminf_{k}\#(D_{k}\bigcap\Omega)/k^{n}\geq(2\pi)^{-n}\int_{\Omega}(i\partial\overline{\partial}\phi)_{n}\label{(I) condition in sampling theorem}\end{equation}
 for any smooth domain $\Omega$ contained in $X(0).$
\end{thm}
\begin{proof}
Given theorem \ref{(I) thm: N is} (applied to the characteristic
function for a set $\Omega$) the proof can be given word by word
as in \cite{li}, \cite{bern}. For completeness we give the argument.
Suppose that the sequence $D_{k}$ is sampling and consider a set
$\Omega$ in $X(0).$ As was explained as a comment to theorem \ref{(I) thm: N is},
the theorem shows that there is subspace of dimension \[
k^{n}(2\pi)^{-n}\int_{\Omega}(i\partial\overline{\partial}\phi)_{n}+o(k^{n})\]
 consisting of functions satisfying the concentration property \begin{equation}
\left\Vert \alpha\right\Vert _{\Omega}^{2}\geq(1-\varepsilon)\left\Vert \alpha\right\Vert _{X}^{2}\label{(I) conc property in pf}\end{equation}
 (we fix some small $\epsilon$).Now the claim is that for any such
concentrated $\alpha$ the sampling property of the sequence $D_{k}$
yields \begin{equation}
\left\Vert \alpha\right\Vert _{_{X}}^{2}\leq Ak^{-n}\sum_{D_{k}\bigcap\Omega_{k}}\left|\alpha(x)\right|^{2},\label{(I) upper estimate on concentrad alpha}\end{equation}
 where $\Omega_{k}$ consist of all points with distance smaller than
$1/\sqrt{k}$ to $\Omega.$ Accepting this for a moment it is easy
to see how the theorem follows. First note that it is enough to prove
the theorem with $\Omega$ replaced with the larger set $\Omega_{k},$
since the number of points in $D_{k}\bigcap(\Omega_{k}-\Omega)$ is
of the order $o(k^{n}).$ To get a contradiction we now assume that
the number of points in $D_{k}\bigcap\Omega_{k}$ is strictly less
than $k^{n}(2\pi)^{-n}\int_{\Omega}(i\partial\overline{\partial}\phi)_{n}+o(k^{n})$
(meaning that condition \ref{(I) condition in sampling theorem} in
the theoren does not hold). But then we can find a non-trivial element
$\alpha$ concentrated on $\Omega$ and vanishing in all the points
in $D_{k}\bigcap\Omega_{k}.$ Indeed, $\alpha$ can be choosen in
a space of dimension of order $k^{n}(2\pi)^{-n}\int_{\Omega}(i\partial\overline{\partial}\phi)_{n}$
and by assumption there are sufficiently few linear conditions to
find such an $\alpha$ vanishing in all the points in $D_{k}\bigcap\Omega_{k}$.
But then \ref{(I) upper estimate on concentrad alpha} forces $\alpha$
to vanish on all of $X,$ which is a contradiction. 

Finally, we just have to show how \ref{(I) upper estimate on concentrad alpha}
follows. For any point $x$ a simple submean inequality gives as in
the begining of section \ref{(I) Section: Asym for B and T}: \[
k^{n}\left|\alpha(x)\right|^{2}\leq C\left\Vert \alpha\right\Vert _{B_{1/\sqrt{k}}(x)}^{2}.\]
 By the separation property of $D_{k}$ we may thus estimate the sum
over $D_{k}\bigcap\Omega_{k}^{c}$ to obtain \[
k^{-n}\sum_{D_{k}\bigcap\Omega_{k}^{c}}\left|\alpha(x)\right|^{2}\leq C\left\Vert \alpha\right\Vert _{\Omega^{c}}^{2}\leq\varepsilon,\]
 where we have used the concentration property \ref{(I) conc property in pf}
in the last step. Hence, we have proved \ref{(I) upper estimate on concentrad alpha},
which finishes the proof of the theorem.
\end{proof}

\part{Harmonic forms}

The aim of the second part of the paper is to generalize the results
in part one to \emph{$\overline{\partial}-$harmonic} $(0,q)-$ forms
with values in $L^{k}.$ We denote the corresponding spaces by $\mathcal{H}^{q}(X,L^{k}),$
which by Hodge's theorem are isomorphic to the Dolbeault cohomology
groups $H^{q}(X,L^{k}).$ The first part was based on the observation
\ref{(I) upper bound on B} that there is always an asymptotic \emph{upper}
bound on the Bergman kernel function $B_{k}(x)$ of the space of holomorphic
sections with values in $L^{k}.$ Furthermore, if the line bundle
$L$ is such that $X(1)$ is empty, then it was shown, using Demailly's
strong Morse inequalities, that the estimate is actually an asymptotic
\emph{equality} (at least in the sense of $L^{1}-$convergence). 

Let us recall the approach to Demailly's inequalities presented in
\cite{berm}. First one shows that \begin{equation}
\dim_{\C}\mathcal{H}_{\leq\nu_{k}}^{q}(X,L^{k})=k^{n}\frac{(-1)^{q}}{\pi^{n}}\int_{X(q)}(\frac{i}{2}\partial\overline{\partial}\phi)_{n}+o(k^{n}),\label{(II) dim formula for Hq}\end{equation}
 for the space $\mathcal{H}_{\leq\nu_{k}}^{q}(X,L^{k})$ spanned by
all eigenforms of $\Delta_{\overline{\partial}}$ with eigenvalues
bounded by $\nu_{k,}$ where $\nu_{k}=\mu_{k}k$ and $\mu_{k}$ is
a certain sequence tending to zero. We will refer to the elements
of the previous space as \emph{low-energy forms.} The dimension formula
\ref{(II) dim formula for Hq} is deduced from the following pointwise
asymptotics for the corresponding Bergman kernel functions:\begin{equation}
B_{\leq\nu_{k}}^{q}(x)=\frac{k^{n}}{\pi^{n}}1_{X(q)}\left|det_{\omega}(\frac{i}{2}\partial\overline{\partial}\phi)_{x}\right|+o(k^{n})\label{(II) Asym for low-e B}\end{equation}
The method of proof is a generalization of the argument used to prove
the upper bound on $B_{k}(x)$ in section \ref{(I) Section: Asym for B and T}
and will not be repeated here. From \ref{(II) Asym for low-e B} one
immediately gets an upper bound on the Bergman kernel functions for
the space of all harmonic forms: \begin{equation}
B_{k}^{q}(x)\leq\frac{k^{n}}{\pi^{n}}1_{X(q)}\left|det_{\omega}(\frac{i}{2}\partial\overline{\partial}\phi)_{x}\right|+o(k^{n})\label{(II) local Morse}\end{equation}
 (in \cite{berm} these bounds were called \emph{local holomorphic
Morse inequalities}). However, the equality \ref{(II) Asym for low-e B}
captures much more of the asymptotic information of the Dolbeault
complex as can be seen in the following way. First observe that the
complex \begin{equation}
\left(\mathcal{H}_{\mu_{k}}^{*}(X,L^{k}),\overline{\partial}\right),\label{(II) complex for fix mu}\end{equation}
 consisting of all eigenforms of $\Delta_{\overline{\partial}}$ with
eigenvalue $\mu_{k}$, forms a finite dimensional subcomplex of the
Dolbeault complex $\left(\Omega^{*}(X,L^{k}),\overline{\partial}\right).$
Indeed, $\overline{\partial}$ commutes with $\Delta_{\overline{\partial}}.$
Moreover, the complex is exact in positive degrees for non-zero $\mu_{k}.$
The Witten $\overline{\partial}-$complex is now defined as the direct
sum of all the complexes \ref{(II) complex for fix mu} with eigenvalue
$\mu_{k}$ less than $\nu_{k}.$ By the Hodge theorem the inclusion
of the Witten $\overline{\partial}-$complex\[
\left(\mathcal{H}_{\leq\nu_{k}}^{*}(X,L^{k}),\overline{\partial}\right)\hookrightarrow\left(\Omega^{*}(X,L^{k}),\overline{\partial}\right),\]
 is a quasi-isomorphism, i.e. the cohomologies of the complexes are
isomorphic and \ref{(II) dim formula for Hq} gives the dimension
of the components of the Witten $\overline{\partial}-$complex. A
homological argument now yields Demailly's \emph{strong} Morse inequalities
for the truncated Euler characteristics of the Dolbeault complex \cite{d1}. 

The main point of this is approach to Demailly's inequalities is to
first prove the Bergman kernel asymptotics \ref{(II) Asym for low-e B}
- the rest of the proof is more or less as Demailly's original proof,
which in turn was inspired by Witten's analytical approach to the
classical real Morse inequalities \cite{wi}. 

Let us now turn to the study of harmonic $(0,q)-$ forms for a fixed
$q,$ i.e. the space $\mathcal{H}^{q}(X,L^{k}).$ As above this space
can be identified with the cohomology groups at degree $q$ of the
Witten $\overline{\partial}-$complex: \begin{equation}
\begin{array}{ccccccc}
 &  & \overline{\partial} &  & \overline{\partial}\\
...\,\rightarrow & \mathcal{H}_{\leq\nu_{k}}^{q-1}(X,L^{k}) & \rightarrow & \mathcal{H}_{\leq\nu_{k}}^{q}(X,L^{k}) & \rightarrow & \mathcal{H}_{\leq\nu_{k}}^{q+1}(X,L^{k}) & \rightarrow\,...\end{array}\label{(II) Witten complex}\end{equation}
 The natural condition on the line bundle $L$ that generalizes the
condition that $X(1)$ is empty, which was used in the study of holomorphic
sections in part one, is that $X(q-1)$ and $X(q+1)$ both are empty. 

\begin{prop}
\label{(II) Prop: equality for Hq}Suppose that $X(q-1)$ and $X(q+1)$
are empty. Then \[
\dim_{\C}\mathcal{H}^{q}(X,L^{k})=k^{n}\frac{(-1)^{q}}{\pi^{n}}\frac{1}{n!}\int_{X(q)}(\frac{i}{2}\partial\overline{\partial}\phi)^{n}+o(k^{n}).\]

\end{prop}
\begin{proof}
Consider the orthogonal decomposition \[
\mathcal{H}_{\leq\nu_{k}}^{q}(X,L^{k})=\mathcal{H}^{q}(X,L^{k})\oplus\mathcal{H}_{+}^{q}(X,L^{k}),\]
 where $\mathcal{H}_{+}^{q}(X,L^{k})$ denotes the eigenspaces corresponding
to positive eigenvalues. According to the dimension formula \ref{(II) dim formula for Hq}
we just have to show that the dimension of $\mathcal{H}_{+}^{q}(X,L^{k})$
is of order $o(k^{n}).$ Since the operator $\overline{\partial}+\overline{\partial}^{*}$
maps $\mathcal{H}_{+}^{q}(X,L^{k})$ injectively into $\mathcal{H}_{\leq\nu_{k}}^{q+1}(X,L^{k})\oplus\mathcal{H}_{\leq\nu_{k}}^{q-1}(X,L^{k})$
the proposition now follows by applying the dimension formula \ref{(II) dim formula for Hq}
again.
\end{proof}
Combining the previous proposition with the upper bound on the Bergman
kernel function $B_{k}^{q}(x)$ of the space of harmonic $(0,q)-$forms
shows that \begin{equation}
k^{-n}B_{k}^{q}(x)\rightarrow\frac{1}{\pi^{n}}1_{X(q)}(x)\left|det_{\omega}(\frac{i}{2}\partial\overline{\partial}\phi)_{x}\right|\label{(II) asym for B}\end{equation}
in $L^{1}(X,\omega_{n}),$ exactly as in the proof of theorem \ref{(I) thm: asym for K}.
When $q$ is zero $B_{k}^{q}(x)$ is the pointwise norm of the restriction
of the Bergman kernel $K(x,y)$ to the diagonal. For $q$ positive,
$K(x,x)$ is (locally) a matrix and $B_{k}^{q}(x)$ is its trace.
But this means that there is now a larger gap between \ref{(II) asym for B}
and the behaviour of $K(x,y)$ - one would rather like to obtain generalizations
of theorem \ref{(I) thm: asym for K} to the matrix elements of $K(x,y).$
To achive this in an invariant way it will be very convenient to think
of $K(x,y)$ as a bundle valued form on $X\times X.$ It will turn
out that $K(x,y)$ is not only localized on the subset corresponding
to $X(q)$ of the diagonal in $X\times X,$ but all of the contribution
to $K(x,x)$ comes from a special direction. This is already clear
from Demailly´s work. In \cite{berm} this fact was shown by reducing
the problem to a model case in $\C^{n}$ where it can be checked explicitely.
In order to treat $K(x,y)$ as a form and to define the special direction
the following formalism will be useful.

\subsection{The super integral and the dagger $\dagger$ (reversed complex conjugation
)}

Consider first a real $m-$dimensional manifold $X.$ Let $f$ be
a differential form on $X,$ i.e. $f$ is an element of $\Omega^{*}(X,\C).$
Then the super integral of $f$ is defined to be the usual integral
of the top degree form of $f,$ i.e\[
\int_{X}f_{1}+f_{2}+...+f_{m}:=\int_{X}f_{m},\]
 where we have decomposed the form $f$ with respect to the degree
grading of $\Omega^{*}(X,\C).$ It is often convenient to think of
the super integral of a form as a double integral in the following
way. Suppose that we are given a volume element, that we write as
$\omega_{m},$ on $X.$ Fix a point $x$ in $X.$ Then $f(x)$ is
element of the exterior algebra over $x$ and we define an {}``integral''
of $f(x)$ by\begin{equation}
\int_{X^{0\mid m}}f(x):=f_{m}(x)/\omega_{m}(x).\label{(II) def: fermionic int}\end{equation}
Next, for a function $f_{0}$ on $X,$ we let \[
\int_{X^{m\mid0}}f(x):=\int_{X}f(x)\omega_{n,}\]
i.e. the usual integral over $X$ of $f_{0}$ with respect to the
volume form $\omega_{n}.$ Then the super integral of a form can be
written as \[
\int_{X}f=\int_{X^{m\mid0}}\int_{X^{0\mid m}}f(x).\]
A word on the notation: In the mathematics litterature the integral
\ref{(II) def: fermionic int} is called the Berezin integral. In
the physics litterature one often thinks of a differential form as
a {}``super'' function of $m$ commuting ({}``bosonic'') variables
$x_{i}$ and $m$ anti-commuting ({}``fermionic'') variables $dx_{i}.$
Taylor-expanding the function $f(x_{1},...,x_{m},dx_{1},...,dx_{m})$
in the anti-commuting variables yields the usual expression of a differential
form. This has been formalized in the theory of super manifolds, where
the super integral corresponds to the integral of a function over
the super manifold $X^{m\mid m},$ that can be obtained from the manifold
$TX^{*}$ by changing the parity along the fibers \cite{dew},\cite{ca}. 

Now assume that $X$ is an n-dimensional complex manifold with a hermitian
metric $\omega$. Consider a $(0,q)-$form $f$ on $X.$ It is well-known
that the squared norm of $f$ may be written as \[
\left\Vert f\right\Vert _{X}^{2}=c_{n,q}\int_{X}f\wedge\overline{f}\wedge\omega^{n-q}/(n-q)!,\]
 in terms of the usual integral, where $c_{n,q}$ is complex constant
needed to make the right hand side real and positive. Using the super
integral we may write the squared norm of \emph{any} form $f$ in
$\Omega^{0,*}(X,\C)$ as \[
\left\Vert f\right\Vert _{X}^{2}=(\frac{i}{2})^{n}\int_{X}f\wedge f^{\dagger}\wedge e^{\omega'},\]
 where $\omega'=-2i\omega$ and where the dagger $^{\dagger}$ is
the linear operator on $\Omega^{*}(X,\C)$ that coincides with the
usual complex conjugation on $\Omega^{1}(X,\C)$ and satisfying \[
\left(\alpha\wedge\beta\right)^{\dagger}=\beta^{\dagger}\wedge\alpha^{\dagger},\]
for any elements $\alpha$ and $\beta$ in $\Omega^{*}(X,\C).$ In
particular $\omega'^{\dagger}=\omega'$ and if if $I=(i_{1},...,i_{q}),$
then \[
dz^{I}\wedge dz^{I\dagger}=(dz^{i_{1}}\wedge\overline{dz}^{i_{1}})\wedge\cdot\cdot\cdot\wedge(dz^{i_{q}}\wedge\overline{dz}^{i_{q}})\]
If there is also given a hermitian line bundle $L$ over $X,$ then
the squared norm of an element $\alpha$ of $\Omega^{0,*}(X,\C)$
may be written as \[
\left\Vert \alpha\right\Vert _{X}^{2}=(\frac{i}{2})^{n}\int_{X}\alpha\wedge\alpha^{\dagger}\wedge e^{-\phi+\omega'}.\]
Note that we are abusing notation here: the function $e^{-\phi}$
representing the fiber metric on $L$ is only defined locally and
$\alpha\wedge\alpha^{\dagger}$ is a $(q,q)-$ form with values in
$L\otimes\overline{L}$ and may not be canonically integrated. However,
the combination $\alpha\wedge\alpha^{\dagger}\wedge e^{-\phi}$ yields
a well-defined global $(q,q)-$form on $X.$ When $X$ is $\C^{n}$
with its standard Euclidean metric form $\omega$and $L$ is the line
bundle with constant positive curvature $\omega',$ the exponent in
the norm above, is equal to \[
\sum_{i}(-z_{i}\overline{z_{i}}+dz_{i}\wedge\overline{dz_{i}})\]
and from the point of view of super manifolds the corresponding Hilbert
space is the space of super functions on $\C^{n\mid n}$ that are
holomorphic in the even variables and anti-holomorphic in the odd
variables.

\subsection{\label{(II) subsection: direction}The direction form $\chi^{q,q}$ }

Fix a point $x$ in $X(q).$ Using the metric $\omega$ we can identify
the curvature two-form $\partial\overline{\partial}\phi_{x}$ at $x$
with a hermitian endomorphism of the fiber over $x$ of the holomorphic
tangentbundle $TX^{1,0}$ in the usual way \cite{gr-ha}. By the definition
of $X(q),$ $\partial\overline{\partial}\phi$ has precisely $q$
negative eigenvalues at $x$ and we denote the complex subspace spanned
by the corresponding eigenvectors by $V(q)_{x}.$ %
\footnote{In the usual real Morse theory the space $V(q)_{x}$ corresponds to
the linearization of the unstable manifold at a critical point.%
} This defines a subbundle $V(q)$ of $TX^{1,0}$over $X(q).$ Denote
the corresponding inclusion map by $i$ and let $\pi$ be the orthogonal
projection of $TX^{1,0}$ onto $V(q).$ On $X(q)$ we define the \emph{direction
form} $\chi^{q,q}$ by \[
\chi^{q,q}:=\pi^{*}i^{*}\omega'_{q}\]
 and extend it by zero to all of $X.$ Locally, $\chi^{q,q}$ can
be expressed in the following way on $X(q).$ Let $(e_{i})$ be a
local orthormal frame of $TX^{1,0}$ such that $e_{1},...,e_{q}$
is a local frame for $V(q)$ and denote by $e^{i}$ the dual $(1,0)-$
forms. Then 

\begin{equation}
\chi^{q,q}:=e^{I_{0}}\wedge e^{I_{0}\dagger}\label{(II) local expression for special direction}\end{equation}
In the sequel, when working with local frames over $X(q),$ we will
assume that $e_{1},...,e_{q}$ are as above. Given a a form $f$ in
$\Omega^{*}(X,\C)$ let $f_{\chi}$ be the function defined by \[
f_{\chi}(x):=(\frac{i}{2})^{n}\int_{X^{0\mid n}}\chi^{q,q}(x)\wedge f(x)\wedge e^{\omega'(x)},\]
 which is just a compact way of saying that on $X(q)$ $f_{\chi}$
is a sum of $f^{0,0}$and all coefficients $f_{JJ}$ such that $J\bigcap I_{0}=\emptyset,$where
$f_{IJ}$ denotes components of the form $f$ with repect to the local
base elements $e^{I}\wedge e^{J\dagger}.$

One final remark: in the notation of the previous section one can
think of $\chi^{q,q}$ as a cut-off function on the super manifold
$X^{n\mid n}.$

\section{Bergman kernels and Toeplitz operators}

Let $(\psi_{i})$ be an orthormal base for a finite dimensional Hilbert
space $\mathcal{H}^{0,q}$ of $(0,q)-$ forms with values in $L.$
Denote by $\pi_{1}$ and $\pi_{2}$ the projections on the factors
of $X\times X.$ The \emph{Bergman kernel} \emph{form} of the Hilbert
space $\mathcal{H}^{0,q}$ is defined by \[
\K(x,y)=\sum_{i}\psi_{i}(x)\wedge\psi_{i}(x)^{\dagger}\]
 Hence, $\K(x,y)$ is a form on $X\times X$ with values in the pulled
back line bundle $\pi_{1}^{*}(L)\otimes\overline{\pi_{2}^{*}(L)}.$
For a fixed point $y$ we identify $\K_{y}(x):=\K(x,y)$ with a $(0,q)-$form
with values in $L\otimes\Omega^{0,q}(X,\overline{L})_{y}.$ The definition
of $\K$ is made so that $\K$ satifisfies the following reproducing
property: \begin{equation}
\alpha(y)=(\frac{i}{2})^{n}\int_{X}\alpha\wedge\K_{y}^{\dagger}\wedge e^{-\phi+\omega'}\label{(II) repr}\end{equation}
for any element $\alpha$ in $\mathcal{H}^{0,q}.$ The restriction
of $\K$ to the diagonal can be identified with a $(q,q)-$form on
$X$ with values in $L\otimes\overline{L}.$ The \emph{Bergman form}
is defined as $\K(x,x)e^{-\phi(x)},$ i.e. \begin{equation}
\B(x)=\sum_{i}\psi_{i}(x)\wedge\psi_{i}(x)^{\dagger}e^{-\phi(x)}\label{II eq def of b form}\end{equation}
and it is a globally well-defined $(q,q)-$form on $X.$ Note that
the Bergman function $B$ is the trace of $\B^{q},$ i.e. \[
B\omega_{n}=c_{n,q}\B\wedge\omega_{n-q}.\]
For a given form $\alpha$ in $\Omega^{0,q}(X,L)$ and a decomposable
form $\theta$ in $\Omega^{0,q}(X)_{x}$ of unit norm, let $\alpha_{\theta}(x)$
denote the element of $\Omega^{0,0}(X,L)_{x}$ defined as \[
\alpha_{\theta}(x)=\left\langle \alpha,\theta\right\rangle _{x}\]
where the scalar product takes values in $L_{x}.$ We call $\alpha_{\theta}(x)$
the \emph{value of $\alpha$ at the point $x,$ in the direction $\theta.$}
Similarly, let $B_{\theta}(x)$ denote the function obtained by replacing
\ref{II eq def of b form} by the sum of the squared pointwise norms
of $\psi_{i,\theta}(x)$. Then $B_{\theta}(x)$ has the following
useful extremal property: \begin{equation}
B_{\theta}(x)=\sup_{\alpha}\left|\alpha_{\theta}(x)\right|^{2},\label{II extremal b}\end{equation}
where the supremum is taken over all elements $\alpha$ in $\mathcal{H}^{q}$
of unit norm. An element $\alpha$ realizing the supremum will be
refered to as an \emph{extremal form for the space $\mathcal{H}^{0,q}$
at the point $x,$ in the direction $\theta.$} The reproducing formula
\ref{(II) repr} may now be written as \begin{equation}
\alpha_{\theta}(y)=(\alpha,\K_{y,\theta})\label{II repr formula}\end{equation}
and we have the following extremal characterization of the Bergman
kernel (which also gives \ref{II extremal b}).

\begin{lem}
\label{(II) lemma: extremal prop of K}Let $\alpha$ be an an extremal
at the point $x$ in the normalized direction $\theta.$ Then \[
\left|\K_{x,\theta}(y)\right|^{2}=\left|\alpha(y)\right|^{2}B_{\theta}(x)\]

\end{lem}
\begin{proof}
Fix the point $x$ and the form $\theta$ in $\Lambda^{n-q,0}(X)_{x}$
and take frames around $x$ such that $\theta=\overline{e^{I}}.$Then
the pair $(x,\theta)$ determines a functional on $\mathcal{H}^{q}:$
\[
\begin{array}{rcl}
\alpha & \mapsto & \alpha_{I}(x)\end{array}.\]
By the reproducing property \ref{II repr formula}\[
\alpha_{I}(x)=(\frac{i}{2})^{n}\int_{X}\alpha\wedge\K_{x,I}\wedge e^{-\phi+\omega'}=(\alpha,\K_{x,I})\]
 for any element $\alpha$ in $\mathcal{H}^{0,q},$ where $\K_{x,I}:=\sum_{i}\psi_{i}^{I}(x)\overline{\psi_{i}}$
is an element of $\mathcal{H}^{q}.$ In terms of a frame at $y$ we
can write \[
\K_{x,I}(y):=\sum_{J}K_{IJ}(x,y)\overline{e^{J}.}\]
 By the reproducing property \ref{II repr formula} $\K_{x,I}/\left\Vert \K_{x,I}\right\Vert $
is an extremal at the point $x$ in the direction $\overline{e^{I}}.$
This means that if $\alpha$ is another extremal at the point $x$
in the direction $\overline{e^{I}}$ then \[
\left|\alpha(y)\right|^{2}=\left|\K_{x,I}(y)/\left\Vert \K_{x,I}\right\Vert \right|^{2},\]
 The previous equality may be written as\begin{equation}
B_{I}(x,x)\left|\alpha(y)\right|^{2}=\left|\K_{x,I}(y)\right|^{2}e^{-\phi(x)},\label{(II) proof of lemma K}\end{equation}
 since the reproducing property \ref{II repr formula} shows that
$\left\Vert \K_{x,I}\right\Vert ^{2}=K_{II}(x,x)$ and by definition
$K_{II}(x,x)e^{-\phi(x)}:=B_{I}(x).$ This proves the lemma.
\end{proof}
Next, denote by $\Omega^{(0)}(X,\C)$ the commutative subalgebra $\bigoplus_{p}\Omega^{p,p}(X,\C)$
of $\Omega^{*}(X,\C).$ For an element $f$ in $\Omega^{(0)}(X,\C),$
we define $T_{f},$ the so called \emph{super Toeplitz operator with
form symbol $f,$} by \begin{equation}
(T_{f}\alpha)(y)=(\frac{i}{2})^{n}\int_{X}f\wedge\alpha\wedge\K_{y}^{\dagger}\wedge e^{-\phi+\omega'}.\label{(II) def of Tf using K}\end{equation}
Equivalently, \begin{equation}
(T_{f}\alpha,\beta)=(\frac{i}{2})^{n}\int_{X}f\wedge\alpha\wedge\beta^{\dagger}\wedge e^{-\phi+\omega'}\label{(II) def of T as bilinear}\end{equation}
 for all elements $\alpha$ and $\beta$ of $\mathcal{H}^{q}.$ Note
that the operator $T_{f}$ is hermitian if $f$ is a real form with
respect to the real structure on $\bigoplus_{p}\Omega^{p,p}(X,\C)$
defined by $\dagger,$i.e. if $f^{\dagger}=f.$ This means that if
$f$ is real-valued in the usual sense i.e. $f$ is an element of
$\Omega^{*}(X,\R),$ then $i^{m}T_{f}$ is hermitian for some integer
$m.$ In the following we will only consider symbols $f$ that are
real with respect to $\dagger.$ 

When studying asymptotic properties of $L^{k},$ all objects introduced
above will be defined with respect to the line bundle $L^{k}.$

\begin{rem}
The term super Toeplitz operator was used in \cite{bor} in a closely
related context. However, the most natural global setting corresponding
to \cite{bor} is obtained by taking the sequence of Hilbert spaces
to be the spaces $\mathcal{H}^{*,0}(X,L^{k}),$ i.e the direct sum
of all harmonic $(q,0)-$forms with values in $L^{k},$ where $q=0,1,..n$
and where $L$ is a \emph{positive} line bundle. Then $\mathcal{H}^{*,0}(X,L^{k})$
is actually the space of all holomorphic forms with values in $L^{k}.$
In particular, the space $\mathcal{H}^{q,0}(X,L^{k})$ may be written
as $H^{0}(X,L^{k}\otimes E_{q}),$ where $E_{q}$ is a holomorphic
vector bundle, so that the analysis for the corresponding Bergman
kernels is reduced to the situation studied in part 1 (twisting with
a fixed vector bundle has only minor effects on the analysis). 
\end{rem}

\section{Asymptotic results for Bergman kernels and Toeplitz operators.}

The next theorem generalizes the bound \ref{(I) local morse q is 0}
in part 1 to low-energy forms on $X,$ i.e. elements of $\mathcal{H}_{\leq\nu_{k}}^{q}(X,L^{k}).$
It is a refined formulation of the local weak holomorphic Morse inequalities
obtained in \cite{berm}.

\begin{thm}
\label{(II) thm: local morse}Fix a point $x$ in $X$ and a direction
form $\theta$ in $\Lambda^{0,q}(X)_{x}.$ Then the following inequality
holds: 

\[
\limsup_{k}\left(k^{-n}\sup\frac{\left|\alpha_{\theta}(x)\right|^{2}}{\left\Vert \alpha\right\Vert _{B_{R_{k}}}^{2}}\right)\leq\frac{1}{\pi^{n}}\left\langle \chi^{q,q},\theta\wedge\theta^{\dagger}\right\rangle _{x}\left|det_{\omega}(\frac{i}{2}\partial\overline{\partial}\phi)_{x}\right|\]
 where the supremum is taken over all elements $\alpha$ of $\mathcal{H}_{\leq\nu_{k}}^{q}(X,L^{k})$
and $R_{k}=\frac{lnk}{\sqrt{k}}.$
\end{thm}
\begin{proof}
In the statement of the theorem concerning local holomorphic Morse
inequalities in \cite{berm} the global norm $\left\Vert \alpha_{k}\right\Vert _{X}^{2}$
was considered. However, the proof given there actually yields the
stronger statement involving $\left\Vert \alpha_{k}\right\Vert _{B_{R_{k}}}^{2}.$
An outline of the argument is as follows. It is enough to prove the
statement for $\theta=\overline{e^{I}}.$ Consider the restriction
of the normalized form $\alpha_{k}$ to the ball $B_{R_{k}}$ centered
at the point $x.$ Let $\beta^{(k)}(z):=k^{-\frac{1}{2}n}\alpha(k^{-\frac{1}{2}}z)$
and extend it by zero to a form on all of $\C^{n}.$ In \cite{berm}
it was shown that one may assume that the sequence $\beta^{(k)}(z)$
tends to a form $\beta$ weakly in $L^{2}(\C^{n}).$ Moreover, the
sequence convergences uniformly with all derivatives on the the unit
ball. This entails that $\beta$ is harmonic with respect to the fiber
metric $\phi_{0}$ (section \ref{(I) Section: Asym for B and T})
on the trivial line bundle in $\C^{n}.$ Thus, \begin{equation}
\limsup_{k}\left|\beta_{I}^{(k)}(0)\right|^{2}/\left\Vert \beta^{(k)}\right\Vert _{B_{lnk}(0)}^{2}\leq\left|\beta_{I}(0)\right|^{2}/\left\Vert \beta\right\Vert _{\C^{n}}^{2}\label{(II) proof of weak Morse}\end{equation}
 where we have used that $\left\Vert \beta\right\Vert _{\C^{n}}^{2}\leq\limsup\left\Vert \beta^{(k)}\right\Vert _{B_{lnk}(0)}^{2}$
thanks to the weak $L^{2}-$ convergence in $\C^{n}.$ In \cite{berm}
the right hand side was shown to be bounded by $\pi{}^{-n}\left|det_{\omega}(\frac{i}{2}\partial\overline{\partial}\phi)_{x}\right|$
if $I=I_{0}$ and equal to zero otherwise. Since the limit in the
left hand side in \ref{(II) proof of weak Morse} equals the limit
of $k^{-n}\left|\alpha_{k}^{I}(x)\right|^{2}/\left\Vert \alpha_{k}\right\Vert _{B_{R_{k}}(x)}^{2},$
this proves the theorem for $\theta=\overline{e^{I}}.$ 
\end{proof}
In \cite{berm} the asymptotics of $B_{k}(x),$ the trace of the Bergman
kernel form $\B_{k}$ associated with the Hilbert spaces $\mathcal{H}_{\leq\nu_{k}}^{q}(X,L^{k}),$
was deduced from the previous theorem. In fact, the proof given there
actually yields the asympototics of the Bergman kernel form itself.
As in part 1 the convergence holds for the Hilbert spaces $\mathcal{H}^{q}(X,L^{k})$
as well, under special conditions on the curvature of $L.$

\begin{thm}
\label{(II) thm: asym for B}Let $\B_{k}$ be the Bergman $(q,q)-$
form of the Hilbert space $\mathcal{H}_{\leq\nu_{k}}^{q}(X,L^{k}).$
Then

\[
k^{-n}\B_{k}(x)\rightarrow\frac{1}{\pi^{n}}\chi^{q,q}\left|det_{\omega}(\frac{i}{2}\partial\overline{\partial}\phi)_{x}\right|\]
pointwise. If $X(q-1)$ and $X(q+1)$ are empty than the convergence
holds in $L^{1}(X,\omega_{n})$ for $\B_{k}^{q}$ associated to the
Hilbert space $\mathcal{H}^{q}(X,L^{k}).$ 
\end{thm}
\begin{proof}
Using the extremal property \ref{II extremal b} of $\B_{k}$, the
upper bound follows immediately from the previous theorem. In particular
$k^{-n}B_{I}(x)$ tends to zero unless $I=I_{0}$ (using frames as
in section \ref{(II) subsection: direction}). Hence, it is enough
to prove the lower bound for the trace $B_{k}(x)$ of $\B_{k}.$ But
this is containded in the asymptotics \ref{(II) asym for B} proved
in \cite{berm} by constructing a sequence of low-energy forms that
become sufficiently large at the point $x,$ when $k$ tends to infinity.
The corresponding result for $\mathcal{H}^{q}(X,L^{k})$ follows just
as in the proof of theorem \ref{(I) Thm: asym for B} in part 1, now
using proposition \ref{(II) Prop: equality for Hq}. 
\end{proof}
The following corollary is obtained just as in part 1: 

\begin{cor}
\label{(II) cor: Tr T}Let $T_{f}$ be the super Toeplitz operator
on the Hilbert space $\mathcal{H}_{\leq\nu_{k}}^{q}(X,L^{k})$ with
symbol form $f.$ Then\[
\lim_{k}\textrm{k}^{-n}\textrm{Tr}T_{f}=\pi{}^{-n}\int_{X}f\wedge\chi^{q,q}\left|det_{\omega}(\frac{i}{2}\partial\overline{\partial}\phi)_{x}\right|\wedge e^{\omega'}.\]
 The right hand side may also be written as \[
(2\pi)^{-n}\int_{X}f_{\chi}(i\partial\overline{\partial}\phi)_{n}.\]
If $X(q-1)$ and $X(q+1)$ are empty then the corresponding result
holds for the Hilbert space $\mathcal{H}^{q}(X,L^{k}).$
\end{cor}
Now we can give the weak convergence of the Bergman kernel stated
in an invariant way.

\begin{thm}
\label{(II) asym for K general}Let $\K_{k}$ be the Bergman kernel
form of the Hilbert space $\mathcal{H}_{\leq\nu_{k}}^{q}(X,L^{k})$
and suppose that $f$ and $g$ are forms in $\Omega^{(0)}(X,\C).$
Then 

\begin{equation}
k^{-n}\int_{X\times X}f(x)\wedge g(y)\wedge\K_{k}(x,y)\wedge\K(x,y)_{k}^{\dagger}\wedge e^{\Phi_{k}(x,y)}\rightarrow(2\pi)^{-n}\int_{X}f_{\chi}g_{\chi}(i\partial\overline{\partial}\phi)_{n},\label{(II) convergence in theorem K general }\end{equation}
 where $\Phi_{k}(x,y)=-k\phi(x)-k\phi(y)+\omega'(x)+\omega'(y).$
If $X(q-1)$ and $X(q+1)$ are empty then the corresponding result
holds for the Hilbert spaces $\mathcal{H}^{q}(X,L^{k}).$
\end{thm}
\begin{proof}
Let us first assume that the Hilbert space is $\mathcal{H}_{\leq\nu_{k}}^{q}(X,L^{k}).$
Consider a subset $U\times V$ of $X,$ where $U$ and $V$ are open
sets in $X(q)$ with associated local frames $e_{U}^{I}$ and $e_{V}^{J}.$
By using a partition of unity it is enough to prove the convergence
with $X\times X$ replaced by $U\times V$ for any such product. Moreover,
by linearity we may assume that $f(x)=F(x)e^{L,L}(x)$ and $g(y)=G(y)e^{M,M}(y),$
where $F$ and $G$ are functions and $e^{I,I}$ is an abrevations
for $e^{I}\wedge e^{I\dagger}.$ Note that $f_{\chi}$ is equal to
$F$ if $L\bigcap I_{0}$ is empty and vanishes otherwise. Using the
special form of $f$ and $g,$ the integral \ref{(II) convergence in theorem K general }
can be written as \begin{equation}
\sum k^{-n}\int_{U\times V}F(x)G(y)\left|K_{IJ}(x,y)\right|^{2}e^{-k\phi(x)-k\phi(y)}\omega_{n}(x)\wedge\omega_{n}(y),\label{(II) aa}\end{equation}
 where the sum is over all $(I,J)$ such that $L\bigcap I$ and $M\bigcap J$
are empty. Let us now show that \begin{equation}
\lim_{k}k^{-n}\int_{U\times V}\left|K_{IJ}(x,y)\right|^{2}e^{-k\phi(x)-k\phi(y)}\omega_{n}(x)\wedge\omega_{n}(y)=0,\label{(II) proof of asym for K a}\end{equation}
unless $U\times V$ is contained in $X(q)\times X(q)$ and $(I,J)=(I_{0},J_{0})$
for the special indices related to the direction form $\chi^{q,q}$
as in \ref{(II) local expression for special direction}. To see this,
assume for example that $U$ is in the complement of $X(q)$ or $I\neq I_{0}$
on $U.$ The integral above is trivially estimated by the limit of
\[
k^{-n}\sum_{J}\int_{U\times X}\left|K_{IJ}(x,y)\right|^{2}e^{-k\phi(x)-k\phi(y)}\omega_{n}(x)\wedge\omega_{n}(y),\]
 Note that the latter integral may be written as the integral of $\left\Vert \K_{x,J}\right\Vert ^{2}e^{-k\phi(x)}$
over all $x$ in $U.$ Now, by the reproducing property \ref{II repr formula}
this integral equals $\int_{U}B_{I}\omega_{n},$ which vanishes if
$U$ is in the complement of $X(q)$ or $I\neq I_{0},$ by theorem
\ref{(II) thm: asym for B}. This proves \ref{(II) proof of asym for K a}.
Using \ref{(II) proof of asym for K a} we may now write the limit
of \ref{(II) aa} as \[
\lim_{k}k^{-n}\int_{U\times V}f_{\chi}(x)g_{\chi}(y)\left|K_{I_{0}J_{0}}(x,y)\right|^{2}e^{-k\phi(x)-k\phi(y)}\omega_{n}(x)\wedge\omega_{n}(y).\]
 Hence to prove the theorem it is is enough to show that if $U\times V$
is contained in $X(q)\times X(q),$ the following holds:\begin{equation}
\lim_{k}k^{-n}\sum_{J}\int_{U\times V}h(x,y)\left|K_{I_{0}J}(x,y)\right|^{2}=\int_{U\bigcap X(q)}h(x,x)(\partial\overline{\partial}\phi)_{n},\label{(II) proof of theorem K b}\end{equation}
for any test function $h(x,y),$ integrating with respect to $e^{-k\phi(x)-k\phi(y)}\omega_{n}(x)\wedge\omega_{n}(y)$
in the left hand side (using \ref{(II) proof of asym for K a} again).
To this end, recall the relation between $K$ and an extremal $\alpha$
at the point $x$ in the direction $\overline{e^{I_{0}}:}$ \[
\sum_{J}\left|K_{I_{0}J}(x,y)\right|^{2}e^{-k\phi(x)-k\phi(y)}=\left|\alpha(y)\right|^{2}e^{-k\phi(y)}B_{I_{0}}(x),\]
given in lemma \ref{(II) lemma: extremal prop of K}. Now the proof
of \ref{(II) proof of theorem K b}, just as the proof of theorem
\ref{(I) thm: asym for K}, is based on the observation that a sequence
of extremals $\alpha_{k}$ at the point $x$ in the direction $I_{0}$
satisfies the localization property \begin{equation}
\lim_{k}\left\Vert \alpha_{k}\right\Vert _{B_{R_{k}}}^{2}=1.\label{(II) proof of theorem K c}\end{equation}
To show this, note that by theorem \ref{(II) asym for B}

\begin{equation}
\lim_{k}k^{-n}\left|\alpha_{k}^{I_{0}}(x)\right|^{2}e^{-k\phi(x)}=\lim_{k}k^{-n}B_{I_{0}}(x)=\pi{}^{-n}\left|det_{\omega}(\frac{i}{2}\partial\overline{\partial}\phi)_{x}\right|\label{(II) proof of theorem K d}\end{equation}
 and by theorem \ref{(II) thm: local morse} 

\begin{equation}
\limsup k^{-n}\left|\alpha_{k}^{I_{0}}(x)\right|^{2}e^{-k\phi(x)}/\left\Vert \alpha\right\Vert _{B_{R_{k}}(x)}^{2}\leq\pi{}^{-n}\left|det_{\omega}(\frac{i}{2}\partial\overline{\partial}\phi)_{x}\right|.\label{(II) proof of theorem K e}\end{equation}
Now \ref{(II) proof of theorem K c} follows from \ref{(II) proof of theorem K d}
together with \ref{(II) proof of theorem K e} just as in the proof
of theorem \ref{(I) thm: asym for K}. Indeed, for a fixed point $x$
the mass of \[
k^{-n}\sum_{J}\left|K_{I_{0}J}(x,y)\right|^{2}e^{-k\phi(x)-k\phi(y)}\]
 considered as a function of $y$ becomes localized close to $y=x,$when
$k$ rw 

Finally, assume that $X(q-1)$ and $X(q+1)$ are empty and consider
the Hilbert space $\mathcal{H}^{q}(X,L^{k}).$ Using the $L^{1}-$convergence
in \ref{(II) asym for B} one sees that \ref{(II) proof of theorem K c}
holds almost everywhere on $X$ for a subsequence of $(\alpha_{k}).$
As in the proof of theorem \ref{(I) thm: asym for K} this is enough
to prove the convergence of $\K(x,y)$ stated in the theorem.
\end{proof}
To formulate the convergence in terms of the local matrix elements
of the Bergman kernel, let $e^{I}$ be a local frame as in section
\ref{(II) subsection: direction} on the open set $U$ in $X.$ Then
we may write the convergence on $U\times U$ in the following suggestive
way \[
\lim_{k}k^{-n}\left|K_{IJ}(x,y)\right|^{2}e^{-k\phi(x)-k\phi(y)}=\delta(x-y)\delta_{IJ}1_{X(q)}(x)1_{q}(I)\pi{}^{-n}\left|det_{\omega}(\frac{i}{2}\partial\overline{\partial}\phi)_{x}\right|\]
 where $1_{q}(I)=1$ if $I=I_{0}$ and zero otherwise. As in part
1 the convergence may also be formulated in terms of (super) Toeplitz
operators:

\begin{cor}
\label{(II) cor: tr TT}Let $T_{f}$ and $T_{g}$ be the super Toeplitz
operators on the Hilbert space $\mathcal{H}_{\leq\nu_{k}}^{q}(X,L^{k})$
with symbol forms $f$ anf $g,$ respectively. Then\[
\lim_{k}\textrm{k}^{-n}\textrm{Tr}(T_{f}T_{g})=\lim_{k}\textrm{k}^{-n}\textrm{Tr}(T_{f_{\chi}g_{\chi}}).\]
If $X(q-1)$ and $X(q+1)$ are empty then the corresponding result
holds for the Hilbert spaces $\mathcal{H}^{q}(X,L^{k}).$
\end{cor}
\begin{proof}
By definition we have that for any $\alpha$ and $\beta$ in $\mathcal{H}^{0,q}$
\[
(T_{f}\alpha,\beta)=\int_{X}f\wedge\alpha\wedge\beta^{\dagger}e^{-k\phi+\omega'}.\]
 Choosing $\alpha=T_{g}\psi$ and $\beta=\Psi$ and and expressing
$T_{f}$ in terms of the Bergman kernel form \ref{(II) def of Tf using K}
gives \[
(T_{f}T_{g}\Psi,\Psi)=\int_{X\times X}f(x)\wedge g(x)\wedge\Psi(x)\wedge\K_{k}(x,y)\wedge\Psi(y)^{\dagger}\wedge e^{\Phi_{k}(x,y)},\]
 where $\Phi_{k}(x,y)=-k\phi(x)-k\phi(y)+\omega'(x)+\omega'(y).$
Finally, if we let $\Psi$ be an orthonormal base element $\Psi_{i}$
and sum over all $i$ the corollary follows from the previous theorem.
\end{proof}
Finally, the following theorem expresses the asymptotic distribution
of the eigenvalues of a super Toeplitz operator in terms of the symbol
of the operator and the curvature of the line bundle $L.$

\begin{thm}
\label{(II) thm: BG}Let $T_{f}$ be the Toeplitz operator with form
symbol $f$ on the Hilbert space $\mathcal{H}_{\leq\nu_{k}}^{q}(X,L^{k}).$
Let $(\tau_{i})$ be the eigenvalues of $T_{f}$ and denote by $d\xi_{k}$
the spectral measure of $T_{f}$ divided by $k^{n},$ i.e. \[
d\xi_{k}:=k^{-n}\sum_{i}\delta_{\tau_{i}},\]
 where $\delta_{\tau_{i}}$ is the Dirac measure centered at $\tau_{i}.$
Then $d\xi_{k}$ tends, in the weak{*}-topology, to the push forward
of the measure $(2\pi)^{-n}(i\partial\overline{\partial}\phi)_{n}$
under the map $f_{\chi},$ i.e. \[
\lim_{i}k^{-n}\sum_{i}a(\tau_{i})=(2\pi)^{-n}\int a(f_{\chi}(x)(i\partial\overline{\partial}\phi)_{n}\]
 for any mesurable function $a$ on the real line. If $X(q-1)$ and
$X(q+1)$ are empty then the corresponding result holds for the Hilbert
spaces $\mathcal{H}^{q}(X,L^{k}).$
\end{thm}
\begin{proof}
As in part 1 we just have to prove the theorem for $a$ equal to the
characteristic function of a half interval, i.e. for the counting
function of $T_{f}.$ Using a partition of unity and the max-min-principle
we may assume that $f$ is supported in a small open set $U$ and
is of the form $Fe^{I,I},$ where $F$ is a function on $U.$ By a
comparison argument it is enough to prove the theorem for $F$ a characteristic
function $1_{\Omega},$ just as in the proof of theorem \ref{(I) thm: N is}
in part 1. Furthermore, as previously we just have to show that \[
\lim_{k}k^{-n}TrT_{f}^{2}=\lim_{k}k^{-n}TrT_{f}.\]
 To this end, observe that for $f=1_{\Omega}e^{I,I}$ clearly $f_{\chi}^{2}=f_{\chi}.$
Hence, corollary \ref{(II) cor: tr TT} shows that $\lim_{k}TrT_{f}^{2}=\lim_{k}TrT_{f_{\chi}},$
which finally is equal to $\lim_{k}TrT_{f},$ by corollary \ref{(II) cor: Tr T}.
\end{proof}
\begin{acknowledgement}
The author wishes to thank his advisor Bo Berndtsson for many interesting
discussions and useful suggestions. 
\end{acknowledgement}


\begin{thebibliography}{10}
\bibitem{berm}Berman, R: Bergman kernels and local holomorphic Morse inequalities
(arXiv.org/abs/math.CV/0211235), to appear in Math Z.
\bibitem{bern}Berndtsson, Bo: Bergman kernels related to hermitian line bundles
over compact comlex manifolds. Explorations in complex and Riemannian
geometry, 1--17, Contemp. Math, 332, Amer. Math. Soc, Providence,
RI, 2003.
\bibitem{bor}Borthwick, D;Klimek,S;Lesniewski, A; Rinaldi, M: Super Toeplitz operators
and nonperturbative deformation quantization of supermanifolds. Comm.
Math. Phys. 153 (1993) 
\bibitem{bouche}Bouche, T: Asymptotic results for Hermitian line bundles over complex
manifolds: the heat kernel approach. Higher-dimensional complex varieties
(Trento, 1994), 67--81, de Gruyter, Berlin, 1996.
\bibitem{bo-gu}Boutet de Monvel, L-Guillemin, V: The spectral theory of Toeplitz
operators. Annals of Mathematical Studies, 99. Princeton University
Press, Princeton, NJ, University of Tokyo Press, Tokyo, 1981
\bibitem{bo-sj}Boutet de Monvel,L-Sjöstrand, J: Sur la singularite des noyaux de
Bergman et Szegö: Journees: Equations aux derivees partielles de Rennes
(1975), 123--164, Asterisque, No. 34-35. Soc. Math. France, Paris,
1976
\bibitem{ca}Cartier, P: DeWitt-Morette, C; Ihl, M; Sämann, C: Supermanifolds-applications
to supersymmetry.
\bibitem{d1}Demailly, J-P: Champs magnetiques et inegalite de Morse pour la d''-cohomologie.,
Ann Inst Fourier, 355 (1985,185-229)
\bibitem{d2}Demailly, J-P: Holomorphic Morse inequalities. Several complex variables
and complex geometry, Part 2 (Santa Cruz, CA, 1989), 93-114
\bibitem{dew}DeWitt, B Supermanifolds. Second edition. Cambridge Monographs on
Mathematical Physics. Cambridge University Press, Cambridge, 1992
\bibitem{do}Donnelly, H: Spectral theory for tensor products of Hermitian holomorphic
line bundles. Math. Z. 245 (2003) no. 1, 31--35
\bibitem{gr-ha}Griffiths, P; Harris, J: Principles of algebraic geometry. Wiley Classics
Library. John Wiley \& Sons, Inc., New York, 1994.
\bibitem{gu}Guillemin, V: Some classical theorems in spectral theory revisited.
219--259 in Seminar on Singularities of Solutions of Linear Partial
Differential Equations. Annals of Mathematics Studies, 91. Princeton
University Press, Princeton, N.J, 1979. 
\bibitem{hö}Hörmander, L: $L^{2}$ estimates and existence theorems for the $\overline{\partial}-$operator.
Acta Math. 113 1965 89--152.
\bibitem{li}Lindholm, N: Sampling in weighted $L^{p}$ spaces of entire functions
in $\C^{n}$ and estimates of the Bergman kernel. J. Funct. Anal.
182 (2001), no. 2, 390--426
\bibitem{ti}Tian, G: On a set of polarized Kähler metrics on algebraic manifolds.
J. Differential Geom. 32 (1990), no. 1, 99--130
\bibitem{wi}Witten, E: Supersymmetry and Morse theory. J. Differential Geom. 17
(1982), no. 4, 661--692. 
\bibitem{ze}Zelditch, S: Szegö kernels and a theorem of Tian. Internat. Math.
Res. Notices 1998, no. 6, 317--331. \end{thebibliography}
\end{document}